\documentclass[12 pt]{article}

\usepackage{amssymb} 
\usepackage{amsmath} 
\usepackage{latexsym}
\usepackage{graphicx}

\usepackage{enumitem}

\graphicspath{ {images/} }
\usepackage{mathtools}

\providecommand{\keywords}[1]{\textbf{\textit{Keywords:}} #1}

\usepackage{tikz}
\usetikzlibrary{calc}
\usepackage{pgf,tikz}
\usetikzlibrary{arrows}
\usepackage{amscd}
\usepackage{relsize}

\usepackage{tikz-cd}



\usepackage{ color, amsmath, amssymb, amsfonts, amscd, graphicx,latexsym, hyperref}
\usepackage{caption, subcaption}
\usepackage[pdf]{pstricks}
\usepackage{makeidx}

\usepackage[affil-it]{authblk}

\usepackage{amsthm}

\newtheorem{theorem}{Theorem}

\newtheorem{proposition}{Proposition}

\newtheorem*{uremark}{Remark}
\newtheorem{lemma}{Lemma}

\newtheorem{definition}{Definition}
\newtheorem*{udefinition}{Definition}
\newtheorem{corollary}{Corollary}

\newcommand{\R}{\mathbb R}

\newcommand{\T}{\mathbb T}

\title{Ideal Triangulation and Disk Unfolding  of a Singular Flat Surface  }
\author{ \.{I}smail Sa\u{g}lam \thanks{Electronic address: \texttt{isaglamtrfr@gmail.com}  }}
\affil{ Adana Alparslan Turkes Science and Technology University}

\date{}

\begin{document}

\maketitle

\begin{abstract}
An ideal triangulation of a singular flat surface is a geodesic triangulation such that its vertex set is equal to the set of singular points of the surface. Using the fact that each pair of points in a surface has a finite number of geodesics having length $\leq L$ connecting them, where $L$ is any positive number, we prove that each singular flat surface has an ideal triangulation provided that the surface has singular points when it has no boundary components, or each of its boundary components has a singular point. Also, we prove that such a surface contains a finite number of geodesics which connect its singular points so that when we cut the surface through these arcs we get a flat disk with a non-singular interior.
\end{abstract}
\keywords{flat surface, raw length spectrum, separatrix, ideal triangulation, unfolding, saddle connection, translation surface}


\section{Introduction}

Let $S_{\frak{g},\frak{n}}$ be a closed orientable surface of genus $\frak{g}\geq0$
with $\frak{n}\geq 0$ punctured points. It is well known that $S_{\frak{g},\frak{n}}$ admits a complete finite-area Riemannian metric of constant curvature equal to -1 when the Euler characteristic $\chi(S_{\frak{g},\frak{n}})$ of $S_{\frak{g},\frak{n}}$ is negative. Indeed, in this case there are many Riemannian metrics of constant curvature equal -1 on $S_{\frak{g},\frak{n}}$. The space of equivalence classes of this metrics is called Teichm\"{u}ller space of $S_{\frak{g},\frak{n}}$ and denoted by $\T(S_{\frak{g},\frak{n}})$. Also, the surface $S_{1,0}$ admits a Riemannian metric of the constant curvature equal to 0. 

It is known that  $\T(S_{\frak{g},\frak{n}})$ is homeomorphic to $\R^{6g-6+2n}$
when $\chi(S_{\frak{g},\frak{n}})<0$. There are several ways to see this. For example, one can form a pants decomposition of $S_{\frak{g},\frak{n}}$
and consider the length parameters and the twist parameters. The coordinates obtained this way are called the Fenchel-Nielsen coordinates. See \cite{thurston3d} or  \cite{FM} for a proof of this fact.

If the surface is not compact, that is $\frak{n}>0$, and $\chi(S_{\frak{g},\frak{n}})<0$,  then there is another way to parametrize the Teichm\"{u}ller space. The main ingredient of this parametrization is an ideal triangulation of a given hyperbolic structure. An ideal triangulation of a hyperbolic surface is a triangulation such that each edge is a geodesic and the vertex set of the triangulation is the set of punctures of the surface. Any hyperbolic structure on such a surface is formed by an ideal triangulation and one may parametrize the Teichm\"{u}ller space by using ideal triangulations. See \cite{AthneseTheret} for details. 

Let $\text{Diff}^+(S_{\frak{g},\frak{n}})$ be the group of all orientation preserving diffeomorphisms of $S_{\frak{g},\frak{n}}$ and $\text{Diff}_0(S_{\frak{g},\frak{n}})$ be the subgroup of diffeomorphisms which are isotopic to the identity. $\text{Diff}^+(S_{\frak{g},\frak{n}})/\text{Diff}_0(S_{\frak{g},\frak{n}})$ 
is called the mapping class group and denoted by $\frak{Mod}(S_{\frak{g},\frak{n}})$. 

If one endowes Teichm\"{u}ller space with the Teichm\"{u}ller's metric, then the mapping class group acts on it by isometries. Also, this action is properly discontinuous. The proof of this fact is non-trivial and the following lemma is one of the main ingredients of the proof. Before stating the lemma, let us recall what the raw length spectrum of a hyperbolic metric is. 

Let $X$ be a complete hyperbolic surface which has finite area. Let $c$ be the isotopy class of a simple closed curve which is  not homotopic to a point or a puncture. We call such a simple closed curve essential. Then it is a well-known fact that $c$ contains a unique geodesic. Let $l_X(c)$ be the length of this simple closed geodesic. Let us define the raw length spectrum of $X$ as the set

$$\frak{rls}(X)=\{l_X(c) \}\subset \R_+,$$

\noindent where $c$ ranges over all isotopy classes of simple closed curves which are essential.

\begin{lemma}[\cite{FM}]
	\label{rawlengthhyp}
	Let $X$ be a hyperbolic srface which has finite area. Then the set $\frak{rls}$ is a closed, discrete subset of $\R$. Furthermore, for each $L \in \R$ the set 
	
	$$\{c: c \ \text{is an isotopy class of an essential simple curve in }\ X \ \text{with} \ l_X(c)\leq L \}$$
	
	\noindent is finite.
	\end{lemma}

Before stating the main results of this paper, let us recall what a singular flat surface is. Roughly speaking, a singular flat surface is a surface which is obtained by gluing Euclidean triangles along their edges appropriately. It is known that any compact surface with genus greater than or equal to one admits a singular flat metric. A geodesic arc on a singular flat surface is a continuous map from a closed interval to the surface which is locally length minimizing.

It is well-known that the space of singular flat surfaces having same topology and prescribed angle data can be identified with the Teichm\"{u}ller space of a punctured surface. See \cite{Trohand}. Therefore, singular flat surfaces form a nice family to study Teichm\"{u}ller theory.

 In \cite{biological}, the authors proved a  result about singular flat surfaces which is similar to Lemma \ref{rawlengthhyp}. The proof is analytic and based on Arzelo-Ascoli Theorem. The authors use this result to obtain a generalized flip algorithm to construct Delaunay triangulation and Voronoi diagram on given singular flat surface. Although this result was proven for only surfaces without boundary, the reasoning of the authors is valid for the singular flat surfaces with a boundary component. Here is the 
 that we mentioned above:

 \begin{proposition}[\cite{biological}]
 	\label{biological}
 	Let $S$ be a compact singular flat surface. For any pair of points $p,q \in S$ and for any $L\geq 0$, the number of geodesics arcs of length $\leq L$ joining $p$ and $q$ is finite. 
 \end{proposition}

We now state the main result of the present paper. It is well known that a compact singular flat surface can be triangulated with Euclidean triangles. A proof of this fact may be found in \cite{TroEnseign}. Indeed one can cut any singular flat surface through piecewise geodesic arcs to obtain a disk with no singular points in its interior, where an interior point is called singular if it has an angle not equal to $2\pi$. See Section \ref{SFS}. Then one can show that such a disk can be triangulated and this triangulation gives us a triangulation of the surface that we started with. However, there are several problems with such a triangulation. First of all, there may be a lot of triangles required and we may not keep track of the number of triangles involved. Also, it does not give us any clue to construct a singular flat surface from a disk or a family of triangles.

\begin{figure}
\hspace*{-2cm} 
		\includegraphics[scale=0.5]{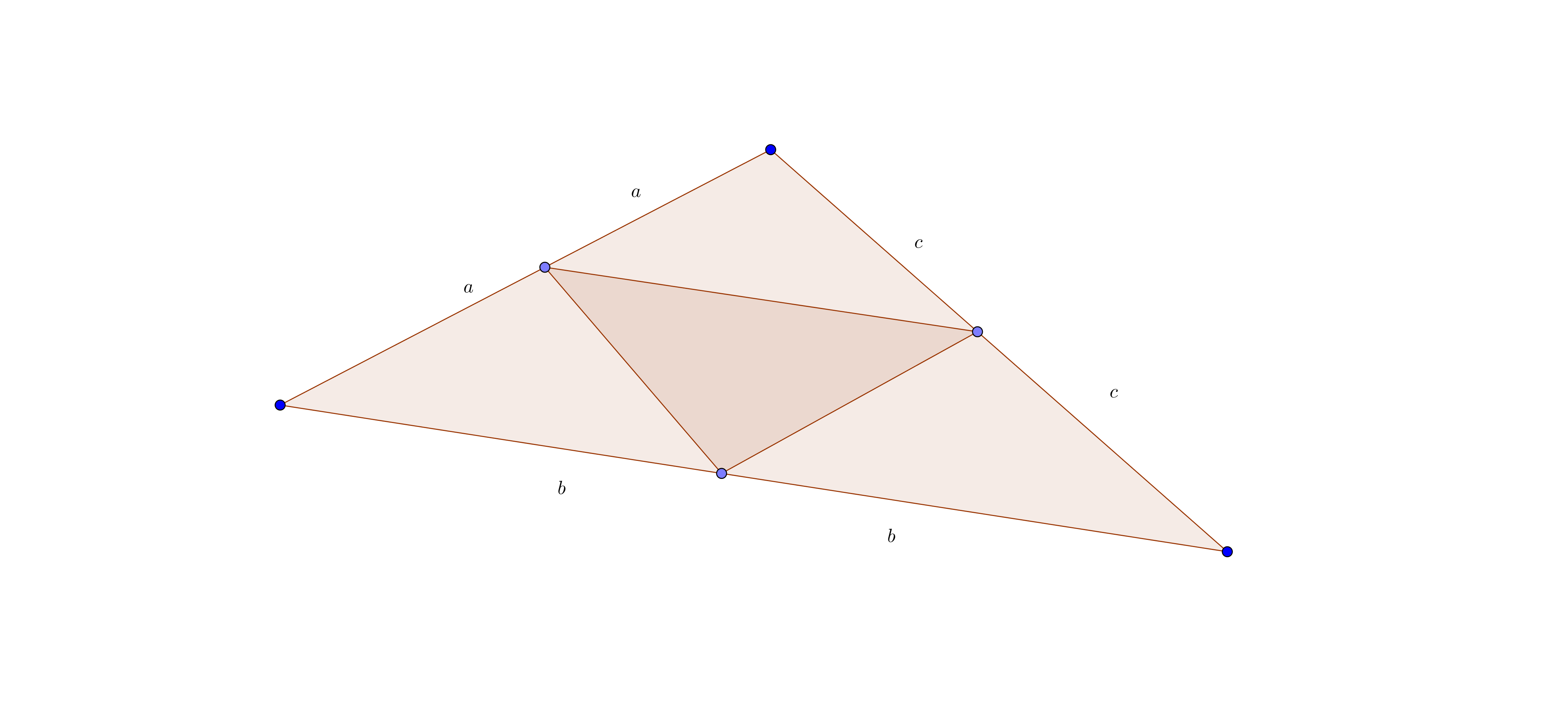}
	
	\caption{A triangulation of the flat sphere with four singular points.}
	\label{tetrahedron}
\end{figure}

If the surface has genus 0 and the angle at each singular point is less than $2\pi$, then it can be triangulated in a nice way. Indeed Thurston \cite{Shapes} proved that such a surface can be triangulated in a way that the vertex set of the triangulation coincides with the set of singular points of the surface. He proved this result by considering Voronoi diagrams in flat surfaces. He used this result to parametrize the space of equivalence classes of singular flat spheres with prescribed angle data. We give an example of such a triangulation in Figure \ref{tetrahedron}. Glue the edges of the triangle with the same label to get a flat sphere with four singular points having angle $\pi$. Then the given triangulation of the triangle induces a triangulation of the tetrahedron.

In this paper, we generalize Thurston's result to arbitrary flat surfaces with or without boundary. More precisely, we prove that any singular flat surface without a boundary component has a triangulation such that the set of singular points coincides with the vertex set of the triangulation provided that the surface has a singular point. If the surface has a boundary component, we show that the surface has such a triangulation
if each of its boundary components has a singular point. Note that we do not use the result of Thurston in the proof. Instead, we use Proposition \ref{biological} to reduce the problem to the case of a disk by an inductive reasoning. Also, we consider non-orientable surfaces with or without boundary.

\begin{figure}
	\hspace*{-2cm} 
	\includegraphics[scale=0.5]{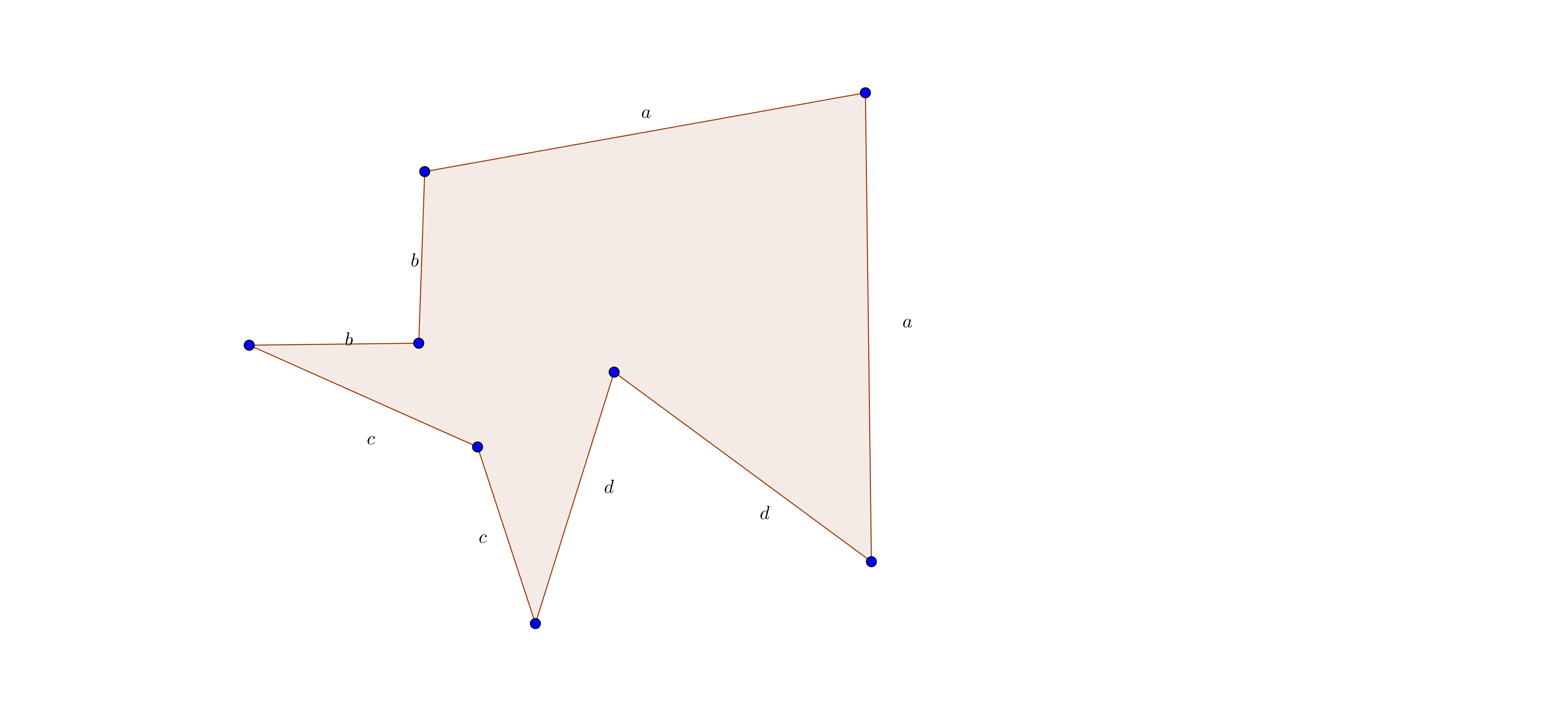}
	
	\caption{A flat sphere from a planer polygon.}
	\label{alexa}
\end{figure}

The second main result of the present paper may be considered as a {\it weak} generalization of the Alexandrov Unfolding Process. See \cite{Alexandrov} and \cite{fillastre}. Alexandrov Unfolding Theorem states that every flat sphere having angle less than $2\pi$ at each singular point can be cut through some geodesic connecting some of its vertices so that the resulted shape is a disk and this disk can be embedded into Euclidean plane. It implies that any flat  sphere with appropriate angle data can be obtained from a planer polygon of a suitable type. We give an example of such a polygon in Figure \ref{alexa}. Glue the edges of the polygon having the same label to obtain a flat sphere with 5 singular points.

We show that any singular flat surface can be obtained from a flat disk with a non-singular interior with a relatively small number of edges. More precisely, we prove  that any singular flat surface  can be cut through a finite number of geodesics which may intersect only at their end points so that the resulting flat surface is a disk provided that the initial singular flat surface has enough singular points.  
 Our result is weak generalization of Alexandrov Unfolding Theorem since it says nothing about embedding of the resulting disk into Euclidean plane. 

\section{Singular Flat Surface}
\label{SFS}
\begin{figure}
	\begin{center}
		\includegraphics[scale=0.6]{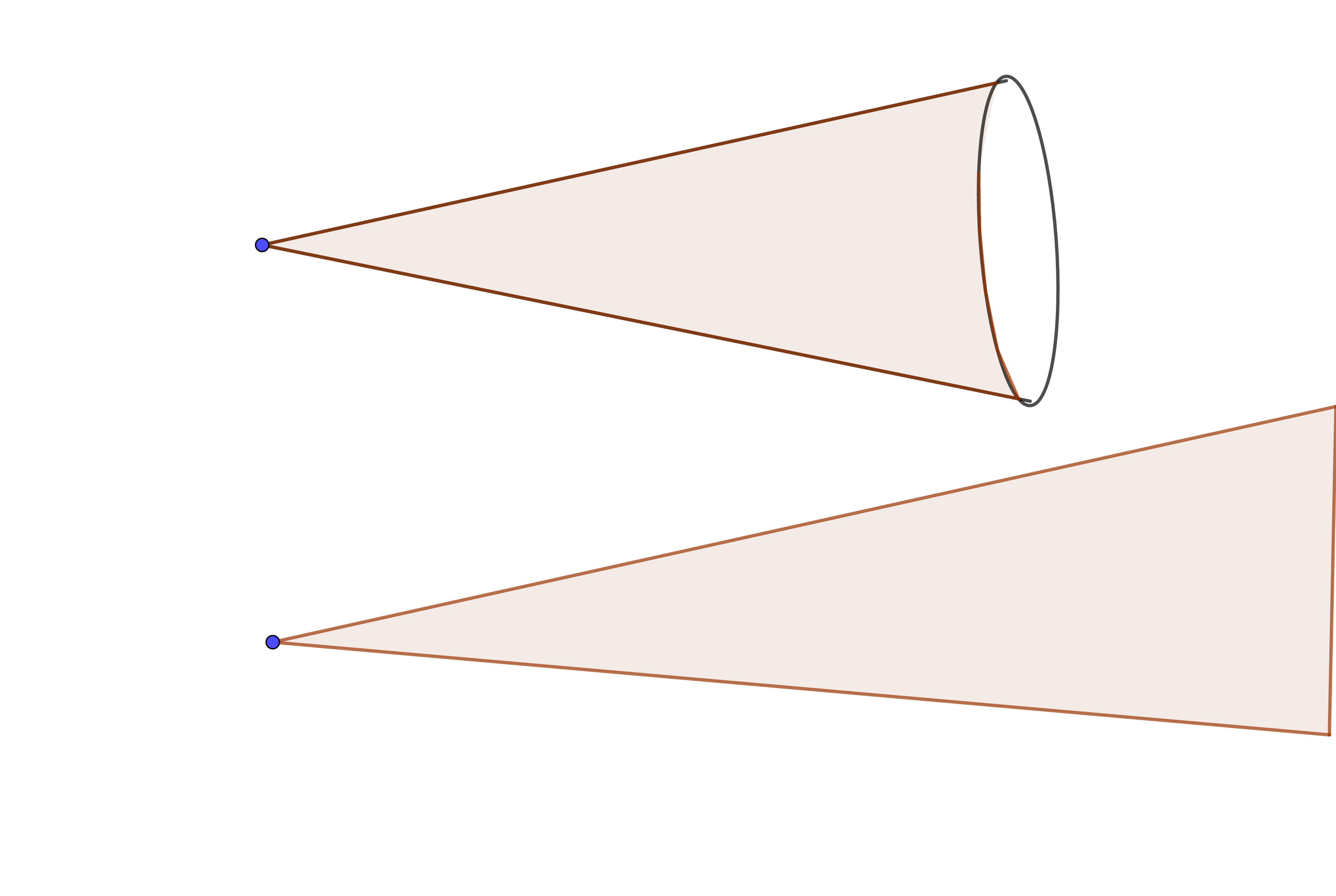}
	\end{center}
	\caption{A cone and a cut of a cone}
	\label{cut}
\end{figure}

Let $S$ be a compact surface perhaps with boundary. A singular flat  metric on $S$ is a singular Riemannian metric so that each interior point of the surface has a neighborhood that is isometric to that of a neighborhood of the apex of  a cone, and each boundary point has a neighborhood isometric to that of a neighborhood of the apex of a cut of a cone. See Figure \ref{cut}. A surface with a singular flat metric is called singular flat surface.

In particular, there is a well defined notion of {\it angle} for each point in a singular flat surface. Let us denote the angle at $p \in S$ by $\theta_p$. An interior point  $p$ of $S$ is called singular  if $\theta_p\neq 2\pi$. Otherwise, it is called non-singular. A boundary point $p$ is called singular if $\theta_p\neq \pi$. Otherwise,  it is called non-singular. The values $2\pi-\theta_p$ and $ \pi-\theta_p$ are called the curvature of $S$ at $p$ if $p$ is an interior point or a boundary point, respectively. Thus a point is non-singular if and only if its curvature is 0.

It is known that a compact singular flat surface can be triangulated. That is, it  can be obtained by gluing Euclidean triangles along their edges appropriately. See \cite{TroEnseign} for a proof. Note that the triangulation given in \cite{TroEnseign} is not necessarily ideal. That is, it may not be true that the set of vertices of the triangulation coincides with the set of singular points of the surface. In \cite{Saglam}, it was proven that any complete singular flat metric in a possibly non-compact surface can be triangulated by finitely many isometry types of Euclidean triangles.

The following formula is called Gauss-Bonnet formula. It's proof may be found in \cite{TroEnseign}.

\begin{proposition}[Gauss-Bonnet Formula]
	Let $S$ be a compact orientable singular flat surface. Then
	
	$$\sum_{p \in S}\kappa_p= 2\pi \chi(S).$$
	
	\end{proposition}

A singular flat surface comes with additional structures. First of all, each path in $S$ has a length. This makes a singular Euclidean surface a length space. That is, we can define a metric on $S$ by considering the infimum of the lengths of the curves that join two points on the surface.  When $S$ is compact, it follows from the general properties of length spaces that for  any two points there is a length minimizing path between them and each homotopy class of paths joining two points has a length minimizing representative. Also, there is a length minimizing geodesic in any free homotopy class of loops. See \cite{AcourseMetric} Note that by a geodesic, we mean a locally length minimizing path.

\section{Labeled Singular Flat Surfaces and Their Raw Length Spectrums }

\label{length-spectrum}

In this section we define labeled singular flat surfaces and their  raw length spectrums.  First, we introduce the notion of the labeled singular flat surface which is essentially a singular flat surface together with a finite set of points which we call the label set. This notion is necessary to prove the existence of a triangulation with desired properties. Proof of this fact is given by induction on the number of singular points and genus  of a given surface, and it is obtained by cutting the surface through arcs. When such a cutting operation 
is performed on a surface, one can get non-singular points from the singular ones. Thus, to keep track of these points, we introduce labeled singular flat surfaces.

\begin{definition}
	A labeled singular flat surface is a compact singular flat surface $S$ together with a finite set of point $P\subset S$ so that
	\begin{itemize}
		\item 
		$P$ is non-empty,
		\item 
		each singular point of $S$ is in $P$,
		\item 
		each boundary component of  $S$ has a non-empty intersection with $P$.
	\end{itemize}
\end{definition}

\begin{definition}
	An arc on a labeled flat surface $(S,P)$ is a geodesic which joins two points of $P$ and whose interior is a subset of $S- P$.
\end{definition}

Note that we assume that constant paths are not arcs. Each arc has a constant speed parametrization and we regard two different constant speed parametrization of the same arc equivalent.

\begin{definition}
	Let $(S, P)$ be a labeled singular flat surface. The raw length spectrum of $(S,P)$ is the set of lengths of arcs on $S$, that is,
	$$\frak{rls}(S,P)=\{l_S(c): c \ \text{is an arc on} \ S \} \subset \R,$$
	
	\noindent where $l_S(c)$ length of the arc $c$.
\end{definition}

\begin{corollary}
	$\frak{rls}(S,P)$ is a discrete and closed subset of $\R$. Indeed, for any $L \geq 0$ the set
	
	$$\{ c: c \ \text{is an arc such that}\  l_S(c)\leq L \}$$
	
	\noindent is finite.
	
	\begin{proof}
		These statements follow immediately from Proposition \ref{biological}.
		\end{proof}
\end{corollary}

Note that when $L$ is sufficiently large $\frak{rls}(S,P)$ contains a positive number. 

\begin{definition}
	Let $(S,p)$ be a labeled flat singular surface. A multiple arc is a path which is a finite union of connected arcs. We define the multiple raw length spectrum as
	
	$$\frak{mrls}=\{l_S(c): c \ \text{is a multiple arc on} \ S\},$$
\noindent	where $l_S(c)$ is the length of $c$ with respect to singular flat metric on $S$.
\end{definition}

The following corollary is immediate.

\begin{corollary}
		$\frak{mrls}(S,P)$ is a discrete and closed subset of $\R$. Indeed, for any $L \geq 0$ the set
	
	$$\{ c: c \ \text{is an multiple arc such that}\  l_S(c)\leq L \}$$
\end{corollary}

Consider the set of geodesics which start and end at same  point of $P$. Call this set $Loop(S,P)$. Note that we do not include 
constant loops to this set. Since each such a geodesic a multiple arc, the following corollary is immediate.

\begin{corollary}\label{loop}
	$Loop(S,P)$ is discrete and and closed. Indeed, for any $L\geq 0$, the set 
	
	$$\{c \in Loop(S,P): \ \text{length of c is less than or equal to} \ L\}$$
	is finite.
\end{corollary}

	A singular flat surface which has a trivial holonomy group is called  a translation surface. A geodesic segment which has non-singular interior and connects two singular points of the surface is called a saddle connection. It is known that the set of holonomy vectors of saddle connections of  translation surface is a discrete subset of $\R^2$. See \cite{handdinamik} for a proof. 

\section{Ideal Triangulation}
In this section, we show that each flat surface with enough singular points have a triangulation whose set of vertices coincides with the set of singular points of the surface. We will make what it means that a surface has enough singular points clear. This amounts to discard some surfaces which have obviously do not have such a triangulation. For example, flat tori can not have such triangulations because they do not have singular points. Also, a singular flat surface with a non-singular boundary component can not have such a triangulation. Indeed, we show that any singular flat surface with boundary having a singular point in each of its boundary component  have such a nice triangulation. Furthermore,  if the surface does not  have boundary components, then it can be triangulated in this manner when it contains singular points.

The proof of existence of such a triangulation is based on induction on the number of boundary components, genus and the number of singular points of the surface. We will cut  a given singular flat surface several times through geodesics that join its singular points to apply the induction hypothesis. It may happen that, if one cuts a singular flat surface through a geodesic starting at a singular point, then one gets non-singular points. This will be the point where the notion of labeled singular flat surface will help us.

Now we define explicitly what we mean by a nice triangulation.

\begin{udefinition}
	Let $S$ be a labeled singular flat surface with the label set $P$. An ideal triangulation of $S$ is a triangulation of $S$ such that
	\begin{enumerate}
		\item 
		each edge of the triangulation is a simple arc,
		\item 
		the vertex set of the triangulation is $P$,
		\item 
	the	interior of each edge does not contain any point of $P$.
	\end{enumerate}
\end{udefinition}
Note that we say that a path $c:[a,b]\to S$ is simple if $t<t'\in [a,b]$ with $c(t)=c(t')$ implies that $t=a$ and $t=b$.
We will prove that any labeled singular flat surface has an ideal triangulation. Note that we do not discard  non-orientable surfaces and surfaces with boundary. 

Some triangles of an ideal triangulation may have two edges and two vertices. For example, one can glue the edges of an isosceles triangle which have equal lengths to get such a triangulation on a singular flat surface. 

\subsection{Orientable surfaces}
Now we start to prove existence of an ideal triangulation. We start with the simplest case; the case of a flat disk. Recall that a flat disk is a labeled 
singular flat surface which is homeomorphic to a closed disk and has a non-singular interior. Note that polygons are examples of flat disks, but it is not true that any flat disk can be embedded into the Euclidean plane. 

\begin{definition}
	A labeled flat disk is a  singular flat disk which is homeomorphic to a closed disk and does not have any labeled points on its interior.
	\end{definition}
\begin{lemma}
	\label{flat disk}
	Each labeled flat disk has an ideal triangulation.
	\begin{proof}
		Any flat disk has at least three singular boundary points. Therefore, if the label set contains three points, then the disk is indeed a triangle, hence it has an obvious ideal triangulation. 
		
		We do induction on the number of labeled points. Assume that a flat disk has a label set which contains  more than three points. Take a point $v$ which is on the boundary and singular. For each other labeled point, consider a length minimizing path joining $v$ and this point. It follows that there exists a labeled point $w$ such that the path joining $v$ and $w$ is not completely contained on the boundary. Note that this is not true if the label set has only three points. Since the path is length minimizing, it follows that it has a subpath such that the end points of this subpath are (not necessarily the same) labeled points and its interior is contained in the interior of the disk. If we cut the disk through this path, then we get two flat disks so that each of them has less labeled points than the original disk. By
		the induction hypothesis, these disks have ideal triangulations, and it is clear that we can glue this disks back to get an ideal triangulation of the original disk.
	\end{proof}
\end{lemma}

Now we consider singular flat surfaces which are homeomorphic to a closed disk possibly with singular interiors.

\begin{lemma}
	\label{flatdisk2}
	A labeled flat surface which is homeomorphic to a closed disk has an ideal triangulation.
	\begin{proof}
		If the disk has no labeled interior points, in Lemma \ref{flat disk}
 we proved that it has an ideal triangulation. Now we do induction on the number of labeled interior points. Since the base case is done, we assume that the disk has labeled interior points. Take a labeled boundary point of the disk and call it $v$. Also take a labeled interior point $w$. There is a (not necessarily unique) length minimizing path between these two points. This path has a subpath which connects a labeled boundary point and a labeled interior point so that interior of it has an empty intersection with $P$. Note that these new points may be different than   $v$ and $w$. Cut the surface through this subpath  and form a new labeled singular flat surface by adding the two new boundary points to the label set. This new labeled singular flat surface has less labeled interior points than the one that we started with. So it has an ideal triangulation by the induction hypothesis. This  triangulation induces an ideal triangulation of the surface that we started with.

	\end{proof}
\end{lemma}

Now we consider the case of the sphere. Note that Thurston \cite{Shapes} also considered this case and obtained an ideal triangulation for each singular flat sphere having non-negative curvature data. His method depends on existence of Veronoi regions on a given flat surface. Our approach is quite different than his approach.

\begin{figure}
	\hspace*{-2cm} 
	\includegraphics[scale=0.5]{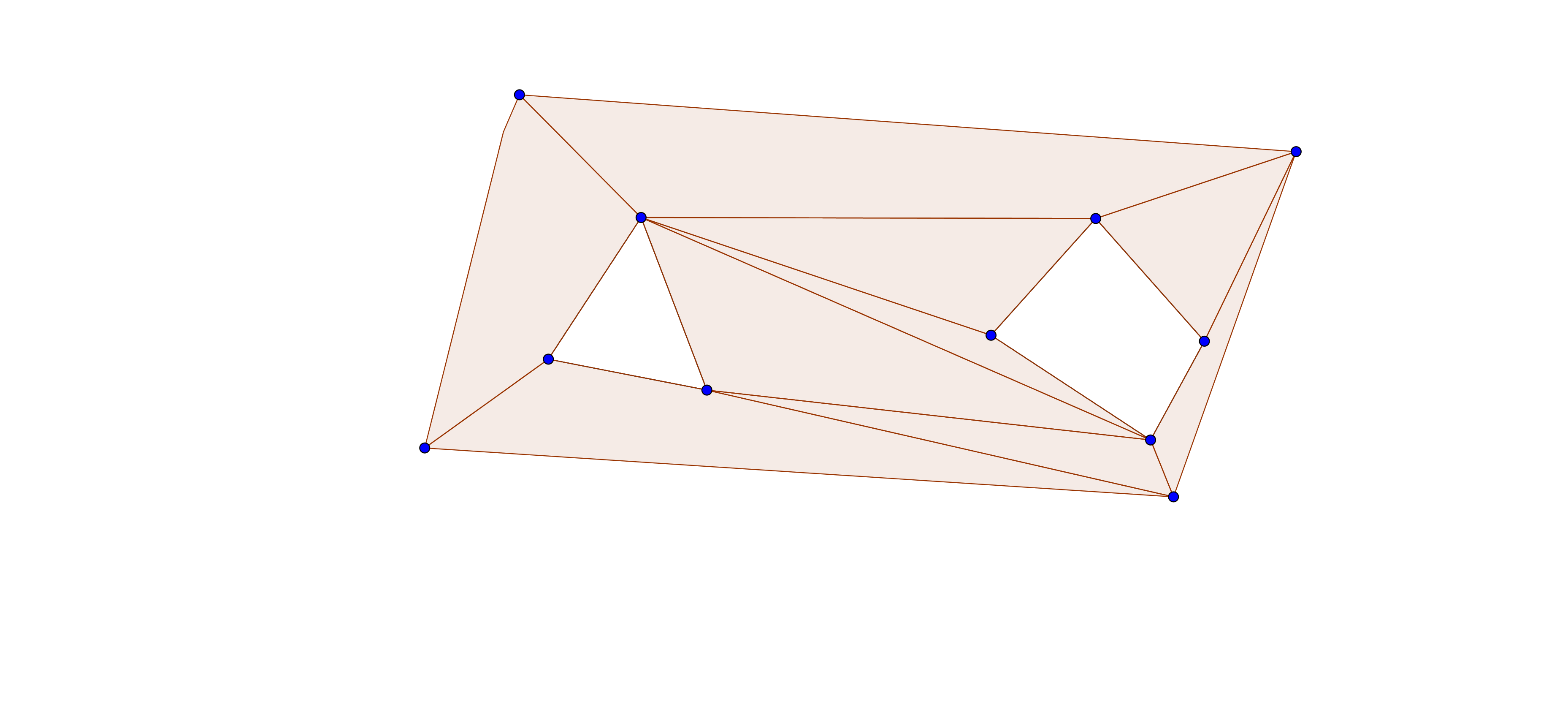}
	
	\caption{A singular flat surface which is homeomorphic to $S^2[3]$.}
	\label{sphere3}
\end{figure}

\begin{lemma}
	\label{sphere}
	Any labeled singular flat surface which is homeomorphic to a sphere has an ideal triangulation.
	
	\begin{proof}
		Note that Gauss-Bonnet formula implies that any singular flat metric on the sphere has at least three singular points. Take two labeled points in such a surface and consider a length minimizing path joining them. Then this path has a subpath which joins two  labeled points and has an interior having non-empty intersection with the label set. Cutting the sphere through this path, we obtain a labeled singular flat surface which is homeomorphic to a closed disk. Lemma \ref{flatdisk2} implies that this disk has an ideal triangulation and it is clear that this triangulation induces an ideal triangulation of the sphere.
	\end{proof}
\end{lemma}
If $S$ is a surface, we denote  the surface obtained by removing interiors of $b$ disjoint closed disks from the interior of $S$ by $S[b]$. Now we consider a singular flat metric on  a surface of genus zero with several boundary components. An example of such a flat surface may be found in Figure \ref{sphere3}.

\begin{lemma}
	\label{sphere-boundary}
	If a labeled singular flat surface is homeomorphic to $S^2[b]$, $b\geq 0$, then it admits an ideal triangulation.
	
	\begin{proof}
		The case $b\leq 1$ was done in Lemma \ref{flatdisk2} and Lemma \ref{sphere}. We do induction on $b$. Assume that $b>1$. Take two distinct boundary components. Let $w$ and $w'$ be labeled points in these boundary components so that they do not lie in the same component. Consider a length minimizing geodesic which connects $w$ and $w'$. Note that if this geodesic makes a non-zero angle between any boundary component at any point distinct than $w$ and $w'$, then this point is a singular point. Furthermore, since the geodesic is length minimizing, it does not intersect itself. Now consider the intersection of $P$ with this geodesic. If we cut the surface through this geodesic, then we get a singular flat surface which is homeomorphic to $S^2[b']$ for some $b'<b$. Each point on the intersection set of the geodesic and $P$ gives us two new points  on the boundary of the resulted surface, we add these points to $P$ to form a new label set. This new surface together with this new label set has an ideal triangulation and this triangulation induces an ideal triangulation of the original labeled flat surface.
	\end{proof} 
\end{lemma}

\begin{figure}
	\hspace*{1cm} 
	\includegraphics[scale=0.5]{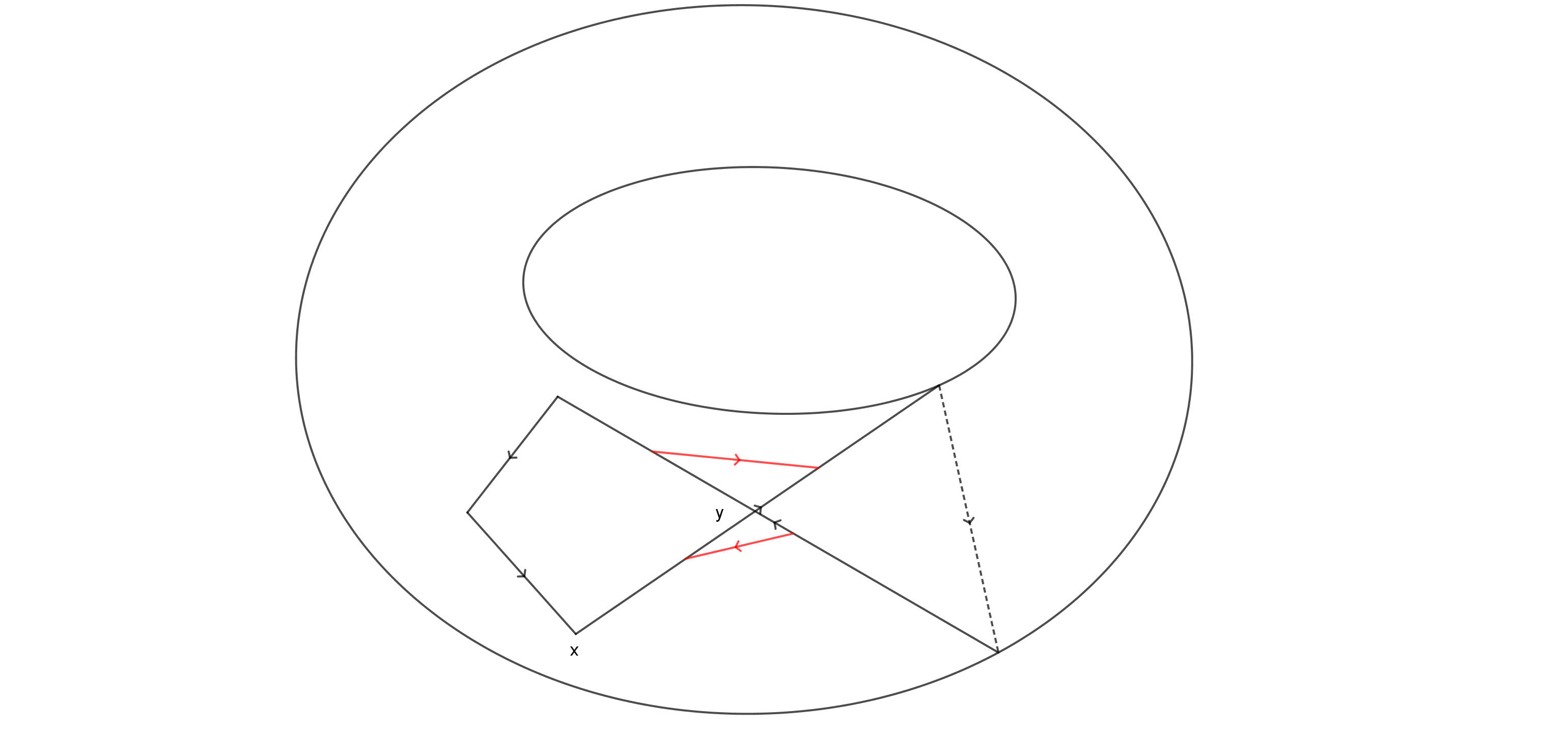}
	
	\caption{The polygonal loop $\alpha$ is described in a torus. Consider the quadrilateral formed by red line segments. If we modify $\alpha$ by replacing its  diagonals by red segments, then we get a path whose length is strictly less than the length of $\alpha$. }
	\label{zor-lemma}
\end{figure}

\begin{lemma}
	\label{simple}
Let $(S,P)$ be a labeled flat surface whose genus is greater than or equal to 1 and which has at most one boundary component.
Let $\frak{S}_P$ be the subset of $Loop(S,P)$ which consists of the curves that are not homotopic to a point and loops winding 
about the boundary component several times. Then the shortest element of $\frak{S}_P$ is a simple geodesic.

\begin{proof}
	Since the genus of the surface is positive, $\frak{S}_P$ is not empty. Note that such an element  is not unique and its existence  follows by  Corollary \ref{loop}. Firstly, assume that the surface has no boundary components. Assume that this geodesic is not simple. Consider the parametric representation $\alpha: [0,1] \to S$ of this geodesic such that $\alpha(0)=\alpha(1)=x \in P$.  Let $b\in (0,1)$ be the smallest value so that there exists $b' >b $ with $\alpha(b)=\alpha(b')$. Call this point $\alpha(b)=\alpha(b')=y$. Then the loop which is obtained by restricting $\alpha$ to the interval $[b,b']$ is not homotopic to a point since the path $\alpha$ is length minimizing. Also the point $y$ can not be in $P$. Indeed, if $y$ is in $P$, then the curve obtaining by restricting $\alpha$ to $[b,b']$ is in $\frak{S}_P$ and has length less than the length of $\alpha$. In particular, $y$ is a non-singular point. Now consider the path obtaining by the juxtaposition of the curves $\alpha\lvert_{[0,b]}$ and $\alpha\lvert_{[b',1]}$. If this loop is homotopic to a point, then there is a piecewise geodesic path  based at $x$ 
	which has length less than  the length of $\alpha$ and which is not homotopically trivial. See Figure \ref{zor-lemma}
	But this means that there is a loop based at $x$ whose length is less then the length of $\alpha$, which is impossible by the choice of $\alpha$.

	Now assume that the surface has one boundary component. Assume that the geodesic is not simple. Consider the parametric representation $\alpha: [0,1] \to S$ of this geodesic such that $\alpha(0)=\alpha(1)=x \in P$. As before, let $b\in (0,1)$ be the smallest value so that there exists $b' >b $ with $\alpha(b)=\alpha(b')$. Call this point $\alpha(b)=\alpha(b')=y$.
    The curve obtaining by restricting $\alpha$ to the interval $[b,b']$ can not be homotopic to a point since $\alpha$ is length minimizing. Assume that it is homotopic to a curve winding about the boundary component several times. Then the curve obtaining by juxtaposing the curves $\alpha\lvert_{[0,b]}$ and $\alpha\lvert_{[b',1]}$ can not be homotopic to a point or the boundary component. But this curve has length than the length of $\alpha$, which is a contradiction. Therefore the curve obtained by restricting $\alpha$ to $[b,b']$ is not homotopic to a point or the boundary. As above it follows that $y$ is not in $P$. In particular, $y$ is not a singular point. Now let $\alpha'$ be the curve joining $y$ and $x$ and which traces the same set with $\alpha\lvert_{[0,b]}$ but in opposite direction. Assume that $\alpha'$ has length less than or equal to $\alpha\lvert_{[b',1]}$. Then the curve obtained by juxtaposing the curves $\alpha\lvert_{[0,b']}$ and $\alpha'$ is not homotopic to the boundary or 
a single point. But it has length less than or equal to the length of $\alpha$. Also, since $y$ is non-singular, we can modify this curve near $y$ to get a curve in same homotopy class but which has strictly smaller length. This contradicts with the choice of  $\alpha$. If the length of 
$\alpha\lvert_{[b',1]}$ is less than the length of $\alpha'$, we can proceed in a similar way to get a contradiction.
\end{proof}
\end{lemma}

\noindent Note that there is a notion of genus for non-orientable surfaces. Lemma \ref{simple} does not discard the case of non-orientable surfaces since in the proof we made no assumption about orientability of $S$.
\begin{proposition}
	\label{orientable}
	Let $(S,P)$ be  an orientable  labeled singular surface. Then $S$ has an ideal triangulation.
	
	\begin{proof}
	Assume that $S$ has more than one boundary components. In that case, as in the proof  Lemma \ref{sphere-boundary}, we may take two points in two different boundary components and cut the surface through a length minimizing geodesic between these  points. Each point in the intersection of $P$ and this geodesic gives us two points in the new surface. By including these points to the label set we can get a new label set. Note that after this operation the number of the boundary components decreases. So we may repeat the process until we get one boundary component. Therefore we assume that $S$ has at most one boundary component. Now we do induction on the genus of $S$. If $S$ has genus zero, then it has an ideal triangulation by Lemma \ref{sphere-boundary}. Assume that genus of $S$ is positive. Lemma \ref{simple} implies that there exists a point $x \in P$ and a simple geodesic loop based at $x$ which is not homotopic to a point or a boundary component. Cut the surface through this geodesic loop. Then there are two possible cases.

\begin{enumerate}
	\item 
	This geodesic is separating. In this case, we obtain two surfaces so that each of them having smaller genus than the original surface. We can form a labeled set for each of these surfaces accordingly. By induction hypothesis, these surfaces admit ideal triangulations. It is clear that these triangulations induce an ideal triangulation of the original surface.
	\item 
	The geodesic is non-separating. In this case we obtain a surface which has genus less than the genus of the original surface. Note that each point on the geodesic induces two points on the boundary of new surface. By adding these points to the labeled set, we may form a new labeled set. By induction hypothesis, this surface has an ideal triangulation and this ideal triangulation induces an ideal triangulation of the original surface.
\end{enumerate}
	\end{proof}
\end{proposition}

The following corollaries are immediate. 

\begin{corollary}
	Let $S$  be a compact orientable surface with boundary together with a singular flat metric. If each boundary component of $S$ has a singular point then $S$ has an ideal triangulation with respect to its set of singular points. 
\end{corollary}

\begin{corollary}
	Let $S$ be a closed orientable surface without boundary and with a singular flat metric. If the set of singular points of $S$ is not empty, then $S$ admits an ideal triangulation with respect to its set of singular points. 
\end{corollary}

\subsection{Non-orientable Surfaces}

Our next objective is to show that any non-orientable labeled flat surface has an ideal triangulation.

\begin{figure}
	\begin{center}
		\includegraphics[scale=0.6]{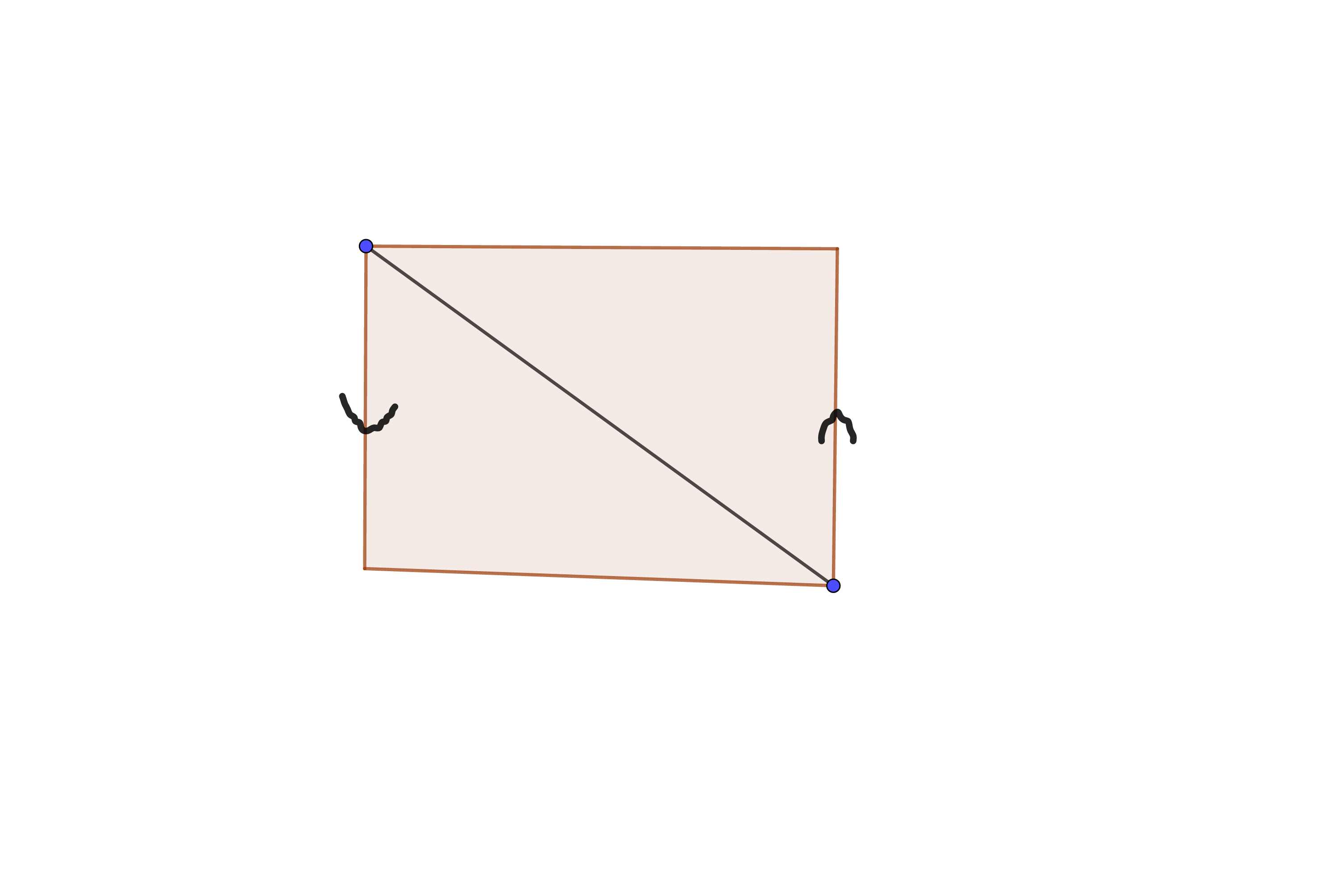}
	\end{center}
	\caption{An ideal triangulation of the M\"{o}bius band with one labeled point.}
	\label{mobius-band}
\end{figure}

\begin{uremark}
	A non-singular M\"{o}bius strip with a label set $P$ having cardinality $n\geq1$ has an ideal triangulation. In Figure \ref{mobius-band}, you can find an example of such a triangulation when the cardinality of $P$ is 1.
\end{uremark}

\begin{lemma}
	\label{ideal-mobius}
	A labeled singular flat surface which is homeomorphic to the M\"{o}bius
	 band admits an ideal triangulation.
	 
	 \begin{proof}
	 Take a labeled boundary point of the surface and consider a length minimizing geodesic which is in the homotopy class of the loop shown in  Figure \ref{mobius-band}. Note that since the geodesic is length minimizing, it does not have self-intersection, that is, it is simple. After this operation we get a disk with a singular flat metric on it. We can form a new label set and get an ideal triangulation of this disk. This triangulation 
	 induces an ideal triangulation of the M\"{o}bius band.
	 \end{proof}
\end{lemma}

\begin{theorem}
	Any compact labeled   singular flat surface admits an ideal triangulation.
	\label{main}
	\begin{proof}
		If the surface is orientable, then the statement was proven in  Proposition \ref{orientable}. Assume that the surface is non-orientable and it is not homeomorphic to the M\"{o}bius band. We do the following operation on the surface:
		
		\begin{enumerate}
			\item 
			Reduce the number of boundary components: if the surface has more than one boundary component, we reduce the number of boundary components to one by cutting the surface through arcs joining its two boundary components repeatedly. 
	\item
	If the surface has less than or equal to one boundary component, then we cut the surface through the loop whose existence is given by Lemma \ref{simple}. 
	
		\end{enumerate}

	After this process, we get one or two singular flat surface. If any of these singular flat surfaces is orientable or homeomorphic to a M\"{o}bius band, then we apply no more cutting 
	operation to it. Otherwise, we apply these operations again and arrange a new set of labeled points for it. Since the surface is of finite type, after finitely many {\it run} of the process we will get a finite number of singular flat surfaces each of which is either homeomorphic to M\"{o}bius band 
	or orientable. By Lemma \ref{ideal-mobius} and Proposition \ref{orientable} these surfaces admit ideal triangulations. These triangulations give an ideal triangulation of the initial surface.
	
	\end{proof}
\end{theorem}

\section{Disk Unfolding}

Consider a singular flat surface which is homeomorphic to the sphere $S^2$ such that curvature at each singular point is positive. Assume that it has $n$ singular points and the curvature at each of these points is positive. Pick a singular point and call it the base point. Consider a length minimizing geodesic for each other singular point which connects the base point and this singular point. Then these geodesics do not have self-intersection and they do not intersect with each other except at the base point. The surface we obtain after cutting the sphere through these geodesics is a flat disk with a non-singular interior. Alexandrov Unfolding Theorem states \cite{Alexandrov} that this flat disk can disk can be embedded isometrically into Euclidean plane. Shortly, it states that any singular flat surface can be obtained from a planer polygon.

In Section \ref{length-spectrum}, we obtained any singular flat surface from a flat disk. But, in the construction, we used an arbitrary number of geodesic segments. Actually, we cut the surface through geodesics which are not necessarily arcs. In this section, we show that we can do this unfolding by using only arcs.

\begin{definition}
	Let $(S,P)$ be a labeled singular flat surface. We say that $(S,P)$ admits a disk unfolding if there exists a finite set of simple arcs of $S$ which can intersect only at their end points so that when we cut the surface through these arcs we get a flat disk with a non-singular interior.
\end{definition}

\begin{lemma}
	\label{sphere-disk-unfolding}
	Any singular flat surface homeomorphic to $S^2[b]$ admits a disk unfolding.
	\begin{proof}
		The proof of this statement is similar to the proof of Lemma \ref{sphere-boundary}. Therefore we omit the proof.
	\end{proof}
\end{lemma}

\begin{proposition}
	\label{orientable-disk-unfolding}
	Any orientable labeled flat surface admits a disk unfolding.
	\begin{proof}
		We do induction on the genus. The case $g=0$ is explained in Lemma \ref{sphere-disk-unfolding}. Assume that the surface has genus greater than zero. If the surface has more than one boundary componet, then we can cut the surface through length minimizing geodesics, which are formed indeed finite unions of arcs, to reduce the number of the boundary components. Therefore we may assume that the surface has at most one boundary component. Lemma \ref{simple} implies that there exists a simple multiple arc which is not homotopic to the boundary component and a point. Now there are two cases to consider.
		
		\begin{enumerate}
			\item 
			The arc is non-separating. If we cut the surface through this arc. We get a surface with smaller genus. By induction hypothesis, this surface admits a disk unfolding and we are done.
			\item 
			The arc is separating. In this case, we obtain two surfaces such that each of them has smaller genus. Each of these surfaces admits disk unfolding. Glue these disks through the corresponding edges, the edges which are induced by the simple arc, to get a disk together with a flat metric. Now this disk admits a disk unfolding and we are done.
		\end{enumerate} 
	\end{proof}
\end{proposition}

\begin{lemma} 
	\label{mobius2}
	Any labeled singular flat surface which is homeomorphic to the M\"{o}bius band admits a disk unfolding. 
	\
\begin{proof}
	The proof Lemma \ref{mobius-band} shows that we can cut this surface through a loop based at its boundary so that the resulting surface is a disk. Then the result follows from Proposition \ref{orientable-disk-unfolding}.
\end{proof}
\end{lemma}

\begin{theorem}
	Any labeled singular flat surface admits a disk unfolding.
	
\begin{proof}
	If the surface is orientable, we proved the statement in Proposition \ref{orientable-disk-unfolding}. Assume that the surface is not orientable.  As in the proof of Theorem \ref{main}, we can cut the surface through the simple arcs 
	to get a finite number of orientable surfaces and M\"{o}bius band. By Proposition \ref{orientable-disk-unfolding} and Lemma \ref{mobius2} each of these surfaces admits a disk unfolding. By gluing these unfoldings accordingly, we may get a disk unfolding.
\end{proof}
\end{theorem}

\end{document}